\documentclass[12pt]{article}
\usepackage{CJK,bbm,amsfonts,mathrsfs,amsthm,amssymb,
amsmath, indentfirst,color}
\usepackage{hyperref}
\def\dy{\displaystyle}
\def\e{\epsilon}
\def\v{\varepsilon}
\def\x{\xi}

\def\a{\alpha}

\def\g{\gamma}
\def\d{\delta}
\def\l{\lambda}

\def\f{\frac}

\def\r{\rho}

\def\di{\displaystyle}
\def\i{\infty}

\def\ef#1 {$(\ref{#1})$}
       \newtheorem{lemma}{\bf Lemma}[section]
       \newtheorem{theorem}[lemma]{\bf Theorem}

       \newtheorem{remark}[lemma]{\bf Remark}

\begin{document}
\begin{CJK}{GB}{gbsn}
\CJKtilde
\title{\bf Zero dissipation limit to rarefaction wave with vacuum for 1-D
compressible Navier-Stokes equations}
\author
 {Feimin Huang,\quad
  Mingjie Li,\quad Yi Wang \\
   {\small \it Institute of Applied Mathematics, AMSS,
   Academia Sinica, Beijing, China.}\\}

\date{}

 \maketitle

\begin{abstract}\noindent It is well-known that  one-dimensional isentropic gas dynamics has two elementary waves, i.e., shock wave and rarefaction wave.
Among the two waves, only the rarefaction wave can be connected to
vacuum. Given a rarefaction wave with one-side vacuum state to the
compressible Euler equations, we can construct a
 sequence of solutions to one-dimensional compressible isentropic
Navier-Stokes equations which converge to the above rarefaction wave
with vacuum as the viscosity tends to zero. Moreover, the uniform
convergence rate is obtained. The proof consists of a scaling
argument and elementary energy analysis based on the underlying
rarefaction wave structures.\end{abstract}

\textbf{Keywords}:~~compressible Navier-Stokes equations, zero
dissipation limit, rarefaction wave, vacuum.

\bigskip

\textbf{    }

 \section{Introduction and main result} \setcounter{equation}{0}
In this paper, we investigate the zero dissipation limit of the
one-dimensional compressible isentropic Navier-Stokes equations
\begin{equation}\label{ns}
\left\{\begin{array}{ll}\di  \rho_{t}+(\rho u)_x=0,&\di x\in \mathbf{R}, t>0,\\[2mm]
\di (\rho u)_t +\big(\rho u^2+p(\rho)\big)_x=\epsilon ~u_{xx},&\\[2mm]
\end{array}\right.
\end{equation}
where $\rho(t,x)\geq0$, $u(t,x)$ and $p$ represent the density, the
velocity and the pressure of the gas, respectively and $\epsilon>0$
is the viscosity coefficient. Here we assume that the viscosity
coefficient $\e$ is a positive constant and the pressure $p$ is
given by the $\g-$law:
$$p(\rho)=\f{\rho^\gamma}{\gamma}$$
 with $\gamma>1$
being the gas constant.

Formally, as $\e$ tends to zero, the limit system of the
compressible Navier-Stokes equations \eqref{ns} is the following
inviscid Euler equations
\begin{equation}\label{eul}
\left\{\begin{array}{l} \rho_{t}+(\r u)_x=0,\\[2mm]
\dy (\r u)_t +(\r u^2+p(\rho))_x=0.
\end{array}\right.
\end{equation}
The Euler system \eqref{eul} is a strictly hyperbolic one for
 $\r>0$ whose characteristic fields are both genuinely nonlinear, that
 is, in the equivalent system
$$
\left(
\begin{array}{l}
\di \r\\
\di u
\end{array}
\right)_t + \left(
\begin{array}{cc}
\di u&\quad \r\\
\di p^\prime(\r)/\r&\quad u
\end{array}\right)\left(
\begin{array}{l}
\di \r\\
\di u
\end{array}\right)_x=0,
$$
the Jacobi matrix
$$
\left(
\begin{array}{cc}
\di u&\quad \r\\
\di p^\prime(\r)/\r&\quad u
\end{array}\right)
$$
has two distinct eigenvalues
$$
\l_1(\r,u)=u-\sqrt{p^\prime(\r)},\qquad
\l_2(\r,u)=u+\sqrt{p^\prime(\r)}
$$
with corresponding right eigenvectors
$$
r_i(\r,u)=(1,(-1)^i\f{\sqrt{p^\prime(\r)}}{\r})^t,\qquad i=1,2,
$$
such that
$$
r_i(\r,u)\cdot \nabla_{\r,u}\l_i(\r,u)=(-1)^i\f{\r
p^{\prime\prime}(\r)+2p^\prime(\r)}{2\r\sqrt{p^\prime(\r)}}\neq
0,\quad i=1,2.
$$
We can define the $i-$Riemann invariant $(i=1,2)$ by
$$
\Sigma_i(\r,u)=u+(-1)^{i+1}\int^\r\f{\sqrt{p^\prime(s)}}{s}ds
$$
such that
$$
\nabla_{(\r,u)}\Sigma_i(\r,u)\cdot r_i(\r,u)\equiv0,\qquad \forall
\r>0,u.
$$

The study of the limiting process of viscous flows when the
viscosity tends to zero, is one of the important problems in the
theory of the compressible fluid.
 When the solution of the inviscid flow is smooth, the zero
dissipation limit  can be solved by classical scaling method.
However, the inviscid compressible flow contains singularities such
as shock and the vacuum in general. Therefore, how to justify the
zero dissipation limit to the Euler equations with basic wave
patterns and/or the vacuum is a natural and difficult problem.

There have been many results on the zero dissipation limit of the
compressible fluid with basic wave patterns without vacuum. For the
system of the hyperbolic conservation laws with artificial viscosity
$$
u_t+f(u)_x=\v u_{xx},
$$
Goodman-Xin \cite{good xin} first verified the viscous limit for
piecewise smooth solutions separated by non-interacting shock waves
using a matched asymptotic expansion method. Later Yu \cite{Yu}
proved it for the corresponding hyperbolic conservation laws with
both shock and initial layers. In 2005, important progress made by
Bianchini-Bressan\cite{B-B} justifies  the vanishing viscosity limit
in BV space even though the problem is still unsolved for the
physical system such as the compressible Navier-Stokes equations.
For the compressible isentropic Navier-Stokes equations \eqref{ns},
Hoff-Liu \cite{hoffliu} first proved the vanishing viscosity limit
for piecewise constant shock even with initial layer. Later Xin
\cite{xin93} obtained the zero dissipation limit for rarefaction
waves without vacuum for both rarefaction wave data and
well-prepared smooth data. Then Wang \cite{Wang-H} generalized the
result of Goodmann-Xin \cite{good xin} to the isentropic
Navier-Stokes equations \eqref{ns}. For the full Navier-Stokes
equations where the conservation of the energy is also involved,
there are also many results on the zero dissipation limit to the
corresponding full Euler system with basic wave patterns without
vacuum. We refer to Jiang-Ni-Sun \cite{jiang} and Xin-Zeng
\cite{Xin-Zeng} for the rarefaction wave, Wang \cite{Wang} for the
shock wave, Ma \cite{Ma} for the contact discontinuity and
Huang-Wang-Yang \cite{Huang-Wang-Yang, Huang-Wang-Yang-1} for the
superposition of two rarefaction waves and a contact discontinuity
 and the superposition of one shock and one rarefaction wave cases.

More recently, Chen-Perepelitsa \cite{chen-P} proved the vanishing
viscosity to the compressible Euler equations for the compressible
Navier-Stokes equations \eqref{ns} by compensated compactness method
for the general case if the far field of the initial values of Euler
system \eqref{eul} has no vacuums. Note that this result is quite
universal since the initial values of the Euler system can contain
vacuum states in the interior domain. Huang-Pan-Wang-Wang-Zhai
\cite{hpwwz} established the corresponding results to the
compressible Navier-Stokes equations \eqref{ns} with
density-dependent viscosity.

Now we turn back to the case of the basic wave patterns with vacuum
states. As pointed out by Liu-Smoller \cite{lius}, among the two
nonlinear waves, i.e., shock and rarefaction waves, to the
one-dimensional compressible isentropic Euler equations \eqref{eul},
only the rarefaction wave can be connected to vacuum. However, to
our knowledge, so far there is no any results on the zero
dissipation limit of the system \eqref{ns} in the case when the
Euler system \eqref{eul} contain the rarefaction wave connected to
the vacuum. In this paper, we investigate this fundamental problem
and want to obtain the decay rate with respect to the viscosity
$\e$. Remark that Perepelitsa \cite{P} consider the time-asymptotic
stability of solutions of 1-d compressible Navier-Stokes equations
\eqref{ns} toward rarefaction waves connected to vacuum in
Lagrangian coordinate and Jiu-Wang-Xin \cite{JWX} study the large
time asymptotic behavior
 toward rarefaction waves for solutions to
  the 1-dimensional compressible
Navier-Stokes equations \eqref{ns} with density-dependent viscosity
for general initial data whose far fields are connected by a
rarefaction wave to the corresponding Euler equations with one end
state being vacuum.

Now we give a description of the rarefaction wave connected to the
vacuum to the compressible Euler equations \eqref{eul}, see also the
references \cite{lius} and \cite{smoller}. For definiteness,
2-rarefaction wave will be considered. If we investigate the
compressible
 Euler system \eqref{eul} with the Riemann initial data
\begin{equation}\label{Riemann}
\left\{
 \begin{array}{rr} \r(0,x)=0,&x<0,\\
(\rho,u)(0,x)=(\r_+,u_+),&x>0,
\end{array}
\right.
\end{equation}
where the left side is the vacuum state and $\r_+>0, u_+$ are
prescribed constants on the right state, then the Riemann problem
\eqref{eul}, \eqref{Riemann} admits a $2-$rarefaction wave connected
to the vacuum on the left side. By the fact that along the
$2-$rarefaction wave curve, $2-$Riemann invariant $\Sigma_2(\r,u)$
is constant in $(x,t)$, we can get the velocity
$u_-=\Sigma_2(\r_+,u_+)$ being the speed of the fluid coming into
the vacuum from the 2-rarefaction wave. This $2-$rarefaction wave
connecting the vacuum $\r=0$ to $(\r_+,u_+)$ is the self-similar
solution $(\r^{r_2},u^{r_2})(\x),~(\x=\f xt)$ of \eqref{eul} defined
by
\begin{equation}
\begin{array}{c}
\di \qquad\qquad\qquad\qquad\quad\r^{r_2}(\x)=0, ~ ~{\rm if}~~\x<\l_2(0,u_-)=u_-,\\
\l_2(\r^{r_2}(\x),u^{r_2}(\x))=\left\{
\begin{array}{ll}
\di
\di \x, &\di {\rm if}~~u_-\leq\x\leq\l_2(\r_+,u_+),\\
\di \l_2(\r_+,u_+), &\di {\rm if}~~\x>\l_2(\r_+,u_+),
\end{array} \right.
\end{array}\label{r2}
\end{equation}
and
\begin{equation}
\Sigma_2(\r^{r_2}(\x),u^{r_2}(\x))=\Sigma_2(0,u_-)=\Sigma_2(\r_+,u_+).\label{r2+}
\end{equation}
Thus we can define the momentum of 2-rarefaction wave by
\begin{equation}m^{r_2}(\x)=\left\{\begin{array}{ll}\r^{r_2}(\x) u^{r_2}(\x),~~~ &{\rm if}~~~\r^{r_2}>0,\\
0,~~~&{\rm if}~~~ \r^{r_2}=0.
\end{array}\right.\label{r2++}
\end{equation}
In the present paper, we want to construct a sequence of solutions
$(\rho^\e,m^\e)(x, t)$ to the compressible Navier-Stokes equations
\eqref{ns} which converge to the 2-rarefaction wave
$(\rho^r_2,m^r_2)(x/t)$ defined above as $\e$ tends to zero. The
effects of initial layers will be ignored by choosing the
well-prepared initial data depending on the viscosity for the
Navier-Stokes equations.

 The main novelty and difficulty of the paper is how to control the
degeneracies caused by the vacuum in the rarefaction wave. To
overcome this difficulty, we first cut off the 2-rarefaction wave
with vacuum along the rarefaction wave curve. More precisely, for
any $\mu>0$ to be determined, the cut-off rarefaction wave will
connect the state $(\r,u)=(\mu,u_\mu)$ and $(\r_+,u_+)$ where
$u_\mu$ can be obtained uniquely by the definition of the
2-rarefaction wave curve. Then an approximate rarefaction wave to
this cut-off rarefaction wave will be constructed through the
Burgers equation. Finally, the desired solution sequences to the
compressible Navier-Stokes equations \eqref{ns} could be established
around the approximate rarefaction wave. The uniform estimates to
the perturbation of the solution sequences around the approximate
rarefaction wave can be got by the following two observations. One
is the fact that the viscosity $\e$ can control the degeneracies
caused by the vacuum in rarefaction waves by choosing suitably
$\mu=\mu(\e)$. In fact, we choose $\mu=\e^{a}|\ln\e|$ with $a$
defined in \eqref{alp} in the present paper. The other observation
is that we can carry out the energy estimates under the a priori
assumption that the perturbation is suitably small in
$H^1(\mathbf{R})$ norm with some decay rate with respect to $\e$ as
$\e$ tends to zero. See \eqref{assump} in the below for the details.
Note that this a priori assumption is natural but is first used in
studying zero dissipation limit to our knowledge. With these two
observations, we can close the a priori assumption and obtain the
desired results.

Now our main result is stated as follows.

\begin{theorem}\label{thm1}
Let $(\rho^{r_2},m^{r_2})(x/t)$ be the 2-rarefaction wave defined by
\eqref{r2}-\eqref{r2++} with one-side vacuum state. Then there
exists a small positive constant $\epsilon_0$ such that for any
$\epsilon\in(0,\epsilon_0)$, we can construct a global smooth
solution $(\rho^\e,m^\e=\rho^\e u^\e)(x,t)$ with initial values
\eqref{inie} to the compressible
Navier-Stokes equation (\ref{ns}) satisfying\\
(1)\begin{equation*}\begin{array}{rl}
(\rho^\e-\rho^{r_2}, m^\e-m^{r_2}),(\rho^\e,m^\e)_x&\dy \in C^0((0,+\infty);L^2(\mathbf{R})),\\
m^\e_{xx}&\dy \in L^2(0,+\infty;  L^2(\mathbf{R})).
\end{array}\end{equation*}
2) As viscosity $\epsilon\rightarrow 0$, $(\rho^\e, m^\e)(x,t)$
converges to $(\rho^{r_2}, m^{r_2})(x/t)$ pointwisely except the
original point $(0,0)$. Furthermore, for any given positive constant
$h$, there exists a constant $C_h>0$, independent of $\epsilon$,
such that
\begin{equation}\label{decay-rate}
\begin{array}{ll}
\dy\sup_{t\geq h}\|
\rho^\e(\cdot,t)-\rho^{r_2}(\frac{\cdot}{t})\|_{L^\infty}\leq
C_h \e^a|\ln\epsilon |,\\
\di \sup_{t\geq h}\|
m^\e(\cdot,t)-m^{r_2}(\frac{\cdot}{t})\|_{L^\infty}\leq
\left\{\begin{array}{ll}\dy C_h\epsilon^{b} |\ln\epsilon |^{-\f12}, \,~~~ &{\rm if}~~~1<\g<3,\\
\dy C_h\epsilon^{\frac{1}{\g+4}} |\ln\epsilon |,~~~&{\rm if}~~~
\g\geq3,
\end{array}\right.
\end{array}
\end{equation}
with the positive constants $a, ~b$ given by
\begin{equation} a=\left\{\begin{array}{ll}\dy\f16 \,~~~ &{\rm if}~~~1<\g\leq2,\\
\dy\frac{1}{\g+4},~~~&{\rm if}~~~ \g>2.
\end{array}\right.\label{aaa}
\end{equation}
and
\begin{equation} b=\left\{\begin{array}{ll}\dy\f18 \,~~~ &{\rm if}~~~1<\g\leq2,\\
\dy\frac{\g+1}{4(\g+4)},~~~&{\rm if}~~~ 2<\g<3.\\[3mm]
\end{array}\right.\label{bbb}
\end{equation}
\end{theorem}
A few remarks are followed.

\begin{remark} Similar result to Theorem~\ref{thm1} is also expected for a one-dimensional
compressible Navier-Stokes equation with density-dependent viscosity
which reads
\begin{equation}\label{sw}
\left\{\begin{array}{l} \rho_{t}+(\rho u)_x=0,\\[2mm]
(\rho u)_t +\big(\rho u^2+\rho^\g\big)_x=\epsilon
\big(\rho^\a u_x\big)_x,\\[2mm]
\end{array}\right.
\end{equation}
with suitable $\a>0$ and $\g>1$. Actually, the system \eqref{sw} can
be derived by Chapman-Enskog expansions from the Boltzmann equation
where the viscosity of the compressible Navier-Stokes equations
depends on the temperature and thus on the density for isentropic
flows. Also, the viscous Saint-Venant system for the shallow water,
derived from the incompressible Navier-Stokes equation with a moving
free surface, is expressed exactly as in \eqref{sw} with $\a=1$ and
$\g=2$,  see $\cite {bn, GP}$.   In this situation, since viscosity
vanishes at vacuum, the convergence rate with respect to $\e$ may
become slower than in Theorem \ref{thm1} and may depend on $\a$ and
$\g$. However, this is left to the forthcoming paper.
\end{remark}

\begin{remark} Our result and method can also be generalized to the 1-D full compressible
 Navier-Stokes equations with the conservation of the energy and the Boltzmann equation with slab symmetry.
 This is left to the forthcoming paper.
\end{remark}

\begin{remark}
It is also interesting to study the zero dissipation limit of
compressible Navier-Stokes equations (\ref{ns}) in the case when the
Euler system (\ref{eul}) has two rarefaction waves with the vacuum
states in the middle. However, it is nontrivial to cut off these
rarefaction waves with vacuum along the corresponding rarefaction
wave curves. In fact, the wave structure containing two rarefaction
waves with the medium vacuum is destroyed and some new wave may
occur in the cut-off precess, which is quite different from the
single rarefaction wave case considered in the present paper.
\end{remark}

\begin{remark}It is noted that in the a priori estimates $(\ref{main})$ below , the estimates for $\phi^2$ from the potential energy holds with the weight
${\bar\r}^{\g-2}$ which is degenerate at vacuum when $\g>2$.
Therefore, the convergence rate obtained in Lemma \ref{len} and thus
in Theorem $\ref{thm1}$ depends on $\g$ when $\g>2$.
\end{remark}

The rest of the paper is organized as follows. In section 2, we
construct a smooth 2-rarefaction wave which approximates the cut-off
rarefaction wave based on the inviscid Burgers equation. And the
proof the Theorem~\ref{thm1} is given in Section 3.

Throughout this paper, $H^l(\mathbf{R}), l= 0,1,2,. . . $,  denotes
the $l$-th order Sobolev space with its norm
$$
\|f\|_l=(\sum^l_{j=0}\|\partial^j_yf\|^2)^\frac{1}{2}, \quad {\rm
and}~\|\cdot\|:=\|\cdot\|_{L^2(dy)},
$$
while $L^2(dz)$ means the $L^2$ integral over $\mathbf{R}$ with
respect to the Lebesgue measure $dz$, and $z=x$ or $y$. For
simplicity, we also write $C$ as generic positive constants which
are independent of time $t$ and viscosity $\e$ unless otherwise
stated.

 \section{Rarefaction waves} \setcounter{equation}{0}
Since there is no exact rarefaction wave profile for the
Navier-Stokes equations \eqref{ns}, the following approximate
rarefaction wave profile satisfying the Euler equations was
motivated by Matsumura-Nishihara \cite{mn86}  and Xin \cite{xin93}.
For the completeness of the presentation, we include its definition
and the properties listed in Lemma \ref{appr}. However, Lemma
\ref{appr} is a little different from \cite{xin93} as stated after
Lemma \ref{appr}.

Consider the Riemann problem for the inviscid Burgers equation:
\begin{align}\label{bur}
\left\{\begin{array}{ll}
w_t+ww_x=0,\\
w(x,0)=\left\{\begin{array}{ll}
w_-,&x<0,\\
w_+,&x>0.
\end{array}
\right.
\end{array}
\right.
\end{align}
If $w_-<w_+$, then the Riemann problem $(\ref {bur})$ admits a
rarefaction wave solution $w^r(x, t) = w^r(\f xt)$ given by
\begin{align}\label{abur}
w^r(\f xt)=\left\{\begin{array}{lr}
w_-,&\f xt\leq w_-,\\
\f xt,&w_-\leq \f xt\leq w_+,\\
 w_+,&\f xt\geq w_+.
\end{array}
\right.
\end{align}
As in \cite{xin93},  the approximate rarefaction wave to the
compressible Navier-Stokes equations \eqref{ns} can be constructed
by the solution of the Burgers equation
\begin{eqnarray}\label{dbur}
\left\{
\begin{array}{l}
\di w_{t}+ww_{x}=0,\\
\di w( 0,x
)=w_\d(x)=w(\f{x}{\d})=\f{w_++w_-}{2}+\f{w_+-w_-}{2}\tanh\f{x}{\d},
\end{array}
\right.\label{(2.11)}
\end{eqnarray}
where $\d>0$ is a small parameter
 to be determined. In fact, we choose $\d=\e^a$ in (\ref{mu}) with $a$ given by \eqref{alp} in the following. Note that the solution $w^r_\d(t,x)$ of the
problem (\ref{(2.11)}) is given by
\begin{equation}\label{b-s}
w^r_\d(t,x)=w_\d(x_0(t,x)),\qquad x=x_0(t,x)+w_\d(x_0(t,x))t.
\end{equation}
And $w^r_\d(t,x)$ has the following properties:

\begin{lemma}
\label{appr} The problem~$(\ref{dbur})$ has a unique smooth global
solution $w_\delta^r(x,t)$ for each $\delta>0$ such that
\begin{itemize}
\item[(1)] $w_-<w_\delta^r(x,t)<w_+, \ \partial_x w^r_\delta(x,t)>0,$
 \ for  $x\in\mathbf{R}, \ t\geq 0,\ \delta>0.$
\item[(2)] The following estimates hold for all $\ t> 0,\ \delta>0$ and
p$\in[1,\infty]$:
\begin{align}
\|\partial_x w^r_\delta(\cdot,t)\|_{L^p}\leq&\dy C
(w_+-w_-)^{1/p}(\delta+t)^{-1+1/p}, \label{w1}
\end{align}
\vspace{-9mm}
\begin{align}\label{w2}
  \|\partial^2_x w^r_\delta(\cdot,t)\|_{L^p}\leq&\dy
C(\delta+t)^{-1}\delta^{-1+1/p},
\end{align}
\vspace{-9mm}
\begin{align}\label{der22}
|\frac{\partial^2 w^r_\delta(x,t)}{\partial
x^2}|\leq\frac{4}{\delta}\frac{\partial w^r_\delta(x,t)}{\partial
x}.
\end{align}
\item[(3)] There exist a constant $\delta_0\in (0,1)$ such that for
$\delta\in(0,\delta_0], t>0$,
$$
\| w^r_\delta(\cdot,t)-w^r(\f\cdot t)\|_{L^\infty}\leq \dy C\delta
t^{-1}\big[\ln (1+t)+|\ln\delta |\big].
$$
\end{itemize}
\end{lemma}

The proof of Lemma \ref{appr} can be found in Xin \cite{xin93}.
However, the description of Lemma \ref{appr} is equivalent to but a
little different from Xin \cite{xin93}. Take the estimation
(\ref{w1}) as an example, which is described by
$$
\|\partial_x w^r_\delta(\cdot,t)\|_{L^p}\leq\dy
C\min\{(w_+-w_-)\d^{-1+1/p}, (w_+-w_-)^{1/p}t^{-1+1/p}\},
\eqno{(\ref{w1})'}
$$
in Xin \cite{xin93}. In fact, two estimations (\ref{w1}) and
$(\ref{w1})'$ are equivalent for fixed wave strength $w_+-w_-$.
However, the advantage of Lemma \ref{appr} is that the energy
estimate can be carried out for all time since there is no
singularity to the approximate rarefaction wave even at $t=0$. While
in Xin's paper \cite{xin93}, the energy estimate must be done in two
time-scalings, that is, finite time and large time, due to the
singularity of the estimations of the approximate rarefaction wave
at $t=0$.

As mentioned in the introduction, we will cut off the 2-rarefaction
wave with vacuum along the wave curve in order to overcome the
difficulty caused by the vacuum,. More precisely, for any $\mu>0$ to
be determined, we can get a state $(\r,u)=(\mu, u_\mu)$ belonging to
the 2-rarefaction wave curve. From the fact that 2-Riemann invariant
$\Sigma_2(\r,u)$ is constant along the 2-rarefaction wave curve,
$u_\mu$ can be computed explicitly by
$u_\mu=\Sigma_2(\r_+,u_+)+\f{2}{\g-1}\mu^{\f{\g-1}{2}}$. Now we get
a new 2-rarefaction wave $(\rho_\mu^{r_2},u_\mu^{r_2}
)(\x),~(\x=x/t)$ connecting the state $(\mu,u_\mu)$ to the state
$(\r_+,u_+)$ which can be expressed explicitly by
\begin{align}\label{rar1}
\lambda_2(\rho_\mu^{r_2},u_\mu^{r_2})(\xi)=\left\{\begin{array}{ll}
\lambda_2(\mu, u_\mu),&\xi< \lambda_2(\mu, u_\mu),\\
\xi,&\lambda_2(\mu, u_\mu)\leq \xi\leq \lambda_2(\rho_+,u_+),\\
\lambda_2(\rho_+,u_+),&\xi>\lambda_2(\rho_+,u_+).
\end{array}
\right.
\end{align}
and
\begin{eqnarray}
\label{rar2}
\Sigma_2(\rho_\mu^{r_2},u_\mu^{r_2})=\Sigma_2(\mu,
u_\mu)=\Sigma_2(\r_+,u_+).
\end{eqnarray}
Correspondingly, we can define the momentum function
$m^{r_2}_\mu=\r^{r_2}_\mu u^{r_2}_\mu.$ It is easy to show that the
cut-off 2-rarefaction wave $(\rho_\mu^{r_2},m_\mu^{r_2} )(x/t)$
converges to the original 2-rarefaction wave with vacuum
$(\rho^{r_2},m^{r_2} )(x/t)$ in sup-norm with the convergence rate
$\mu$ as $\mu$ tends to zero. More precisely, we have

\begin{lemma}\label{cut-off}
There exist a constant $\mu_0\in (0,1)$ such that for
$\mu\in(0,\mu_0], t>0$,
\begin{center}
$\| (\rho_\mu^{r_2},m_\mu^{r_2} )(\cdot/t)-(\rho^{r_2},m^{r_2}
)(\cdot/t)\|_{L^\infty}\leq \dy C\mu.$
\end{center}
\end{lemma}
The proof of Lemma \ref{cut-off} can be obtained directly from the
explicit solution formula of rarefaction waves, so we omit it for
brevity.

Now the approximate rarefaction wave
$(\bar{\r}_{\mu,\d},\bar{u}_{\mu,\d})(x,t)$ of the cut-off
2-rarefaction wave $(\rho_\mu^{r_2},u_\mu^{r_2})(\f xt)$ to
compressible Navier-Stokes equations $(\ref {ns})$ can be defined by
\begin{eqnarray}
\left\{
\begin{array}{l}
\di w_+=\l_2(\r_+,u_+),\quad w_-=\l_2(\mu,u_\mu), \\
\di w_\d^r(t,x)= \l_2(\bar{\r}_{\mu,\d},\bar{u}_{\mu,\d})(t,x),\\
\di
\Sigma_2(\bar{\r}_{\mu,\d},\bar{u}_{\mu,\d})(x,t)=\Sigma_2(\rho_+,u_+)=
\Sigma_2(\mu,u_\mu),
\end{array} \right.\label{au}
\end{eqnarray}
where $w_\delta^r$ is the solution of Burger's equation $(\ref
{dbur})$ defined in \eqref{b-s}. From then on, the subscription of
$(\bar\r_{\d,\mu}, \bar u_{\d,\mu})(x,t)$ will be  omitted as
$(\bar\r,\bar u)(x,t)$ for simplicity. Then the approximate cut-off
2-rarefaction wave $(\bar \r,\bar u)$ defined above satisfies
\begin{equation}\label{ar}
\left\{\begin{array}{ll}\di  \bar\rho_{t}+(\bar\rho\bar u)_x=0,\\[2mm]
\di (\bar\rho \bar u)_t +\big(\bar\rho\bar u^2+p(\bar\rho)\big)_x=0,\\[2mm]
\end{array}\right.
\end{equation}
and the properties of the approximate rarefaction wave $(\bar\r,\bar
u)$ is listed without proof in the following Lemma.

\begin{lemma}
\label{appu} The approximate cut-off 2-rarefaction wave $(\bar
\r,\bar u)$ defined in \eqref{au} satisfies the following
properties:
\begin{itemize}
\item[(i)] $\bar u_x(x,t)=\f2{\g+1}(w_\d^r)_x>0,$
 \ for  $x\in\mathbf{R}, \ t\geq 0;$\\
$\bar\r_x=\bar\r^{\f{3-\g}{2}}\bar u_x$, and
$\bar\r_{xx}=\bar\r^{\f{3-\g}{2}}\bar
u_{xx}+\f{3-\g}{2}\bar\r^{2-\g}(\bar u_x)^2$.
\item[(ii)] The following estimates hold for all $\ t> 0,\ \delta>0$ and
p$\in[1,\infty]$:
\begin{equation*}\begin{array}{l}
\| \bar u_x(\cdot,t)\|_{L^p}\leq
C (w_+-w_-)^{1/p}(\delta+t)^{-1+1/p},\\[4mm]
  \|\bar u_{xx}(\cdot,t)\|_{L^p}\leq
C(\delta+t)^{-1}\delta^{-1+1/p}.
\end{array}\end{equation*}
\item[(iii)] There exist a constant $\delta_0\in (0,1)$ such that for
$\delta\in(0,\delta_0], t>0$,
\begin{equation*}\begin{array}{ll}
\| (\bar\r-\rho^{r_2}_\mu, \bar
u-u^{r_2}_\mu)(\cdot,t)\|_{L^\infty}\leq \dy C\delta t^{-1}\big[\ln
(1+t)+|\ln\delta |\big].
\end{array}\end{equation*}
\end{itemize}
\end{lemma}

\section{Proof of Theorem~\ref{thm1}}\setcounter{equation}{0}

To prove Theorem $\ref{thm1}$, the solution $(\rho^\e, u^\e)$ is
constructed as the perturbation around the approximate rarefaction
wave $(\bar\rho, \bar u)$.  Consider the Cauchy problem for
\eqref{ns} with smooth initial data
\begin{equation}
\begin{aligned}\label{inie}
(\rho^\e, u^\e)(x,t=0)=(\bar\r, \bar u)(x,0).
\end{aligned}
\end{equation}
Then we introduce the perturbation
\begin{equation}\label{pert}
(\phi, \psi)(y,\tau)=(\rho^\e, u^\e)(x,t)- (\bar\r, \bar u)(x,t),
\end{equation}
where $y,\tau$ are the scaled variables as
\begin{equation}\label{tao}
y=\frac{x}{\epsilon},\quad  \tau=\frac{t}{\epsilon},
\end{equation}
and $(\rho^\e, u^\e)$ is assumed to be the solution to the problem
 $(\ref {ns})$. For the simplicity of the notation, the superscription of $(\rho^\e,
 u^\e)$ will be omitted
  as $(\rho, u)$ from now on if there is no confusion of the notation. Substituting $(\ref{pert})$ and $(\ref{tao})$
into $(\ref {ns})$  and using the definition for $(\bar\r, \bar u)$,
we obtain
\begin{align}
\phi_\tau+\r\psi_y+u\phi_y=-f\label{mass},
\\
\r\psi_\tau+\r u\psi_y+p'(\r)\phi_y- \psi_{yy}=-g, \label{mome}
\end{align}
\begin{equation}\label{init}
(\phi, \psi)(y,0)=0,
\end{equation}
where
\begin{equation}\label{g}
\left\{
\begin{array}{l}
f={\bar u}_y\phi+{\bar\r}_y\psi,
 \\[2mm]
\dy g=-\bar u_{yy}+\r\psi \bar u_y+
\bar\rho_y\Big[p'(\r)-\frac{\r}{\bar\r} p'(\bar\r)\Big].
\end{array}\right.
\end{equation}

We seek a global (in time) solution $(\phi,\psi)$ to the
reformulated problem $(\ref{mass})-(\ref{init})$. To this end, the
solution space for $(\ref{mass})-(\ref{init})$ is defined by
\begin{align}
X(0,\tau_1)=\Big\{(\phi,\psi)\Big|&(\phi,\psi)\in
C^{0}([0,\tau_1];H^1(\mathbf{R})),\nonumber\quad
\phi_y\in L^{2}(0,\tau_1;L^2(\mathbf{R})),\nonumber\\
&\psi_y\in L^{2}(0,\tau_1;H^1(\mathbf{R})) \Big\}\nonumber
  \end{align}
with $0<\tau_1\leq+\infty$.

\begin{theorem}\label{thm31}
The problem $(\ref{mass})-(\ref{init})$ admits a unique
global-in-time solution $(\phi,\psi)\in X(0,+\i)$. Furthermore,
there exist positive constants $\e_0$ and $C$ independent of $\e$,
such that if $0<\e\leq\e_0$,~ then
\begin{equation}\label{main-1}\begin{array}{ll}
&\dy\sup_{\tau\in[0,+\i]}\int_{\mathbf{R}} \Big(\bar\r\psi^2+{\bar\r}^{\g-2}\phi^2+\phi^2_y+\psi^2_y\Big)(\tau,y)dy\\[4mm]
+&\dy\int^{+\i}_{0}\int_{\mathbf{R}}
\Big[\psi_y^2+{\bar\r}^{\g-2}\bar u_y\phi^2+\bar\r \bar
u_y\psi^2+\bar \r^{\g-3}\phi^2_y+\frac{\psi^2_{yy}}{\bar\r}\Big]
dyd\tau \leq C\e^{(1/2-a)}|\ln\e|^{-1/2}.
\end{array}\end{equation}
Consequently,
\begin{equation}\label{rate-1}
\begin{array}{ll}
\di \sup_{0\leq\tau\leq+\i}\|\phi(\cdot,\tau)\|_{L^\infty} \leq
\left\{\begin{array}{ll}\dy C\e^{1/6}|\ln\e|^{-1/4},~~~ &{\rm if}~~~1<\g\leq2,\\
\dy C\e^\frac{1}{\g+4}|\ln\e|^{(1-\g)/4},~~~&{\rm if}~~~ \g>2,
\end{array}\right. \\[5mm]
\di \sup_{0\leq\tau\leq+\i}\|\psi(\cdot,\tau)\|_{L^\infty} \leq
\left\{\begin{array}{ll}\dy C\e^{1/8}|\ln\e|^{-1/2},~~&{\rm if}~~~1<\g\leq2,\\
\dy C\e^\frac{\g+1}{4(\g+4)}|\ln\e|^{-1/2},~~&{\rm if}~~~ \g>2.
\end{array}\right.
\end{array}
\end{equation}
\end{theorem}

In what follows, the analysis is always carried out  under the a
priori assumptions
\begin{align}\label{assump}
\sup_{0\leq\tau\leq\tau_1} \|\phi(\cdot,\tau)\|_{L^\infty}\leq \e^{a},~~~
\sup_{\tau\in[0,\tau_1]}\|\psi_y\|\leq 1,
\end{align}
with $a$ given by
\begin{align}\label{alp}
a=\left\{\begin{array}{ll}
\dy\frac{1}{6},& 1< \gamma\leq2,\\
\dy\frac{1}{\gamma+4},&  \gamma>2.
\end{array}
\right.
\end{align}
Take
\begin{equation}\label{mu}
\mu=\e^{a}|\ln\e|, \qquad\d=\e^{a},
\end{equation}
in the sequel. Then it follows that $\mu\geq2 \e^{a}$ if $\e\ll1$.
Under the a priori assumption \eqref{assump}, we can get
\begin{align}\label{rup}
\frac{\bar\r}{2}\leq \r\leq\frac{3\bar\r}{2}.
\end{align}
In fact, if $\e\ll1$, then
\begin{align}
\r=\bar\r+\phi\geq \bar\r-\|\phi\|_{L^\infty}\geq\bar\r-\e^{a}\geq
\bar\r-\frac{1}{2}\mu\geq\frac{\bar\r}{2},\\
\r=\bar\r+\phi\leq \bar\r+\|\phi\|_{L^\infty}\leq\bar\r+\e^{a}\leq
\bar\r+\frac{1}{2}\mu\leq\frac{3\bar\r}{2}.
\end{align}
Moreover, under the a priori assumption \eqref{assump}, it holds
that
\begin{align}\label{eqr}
C_1{\bar\r}^{\g-2}\phi^2\leq p(\r)-p(\bar\r)-p'(\bar\r)\phi\leq
C_2{\bar\r}^{\g-2}\phi^2.
\end{align}
where $C_1, C_2$ are positive constants independent of $\e$.

Since the proof for the local existence of the solution to
$(\ref{mass})-(\ref{init})$ is standard, we omit it for brevity. To
prove Theorem~\ref{thm31}, it is sufficient to prove the following a
priori estimates.

\begin{lemma}\label{len}
\textbf{ (A priori estimates)}\ \ Let $\gamma>1$ and $(\phi,\psi)\in
X(0,\tau_1)$ be a solution to the problem
$(\ref{mass})-(\ref{init})$. Then under the a priori assumption
\eqref{assump}, there exist positive constants $\e_0$ and $C$
independent of $\e$, such that if $0<\e\leq\e_0$,~ then
\begin{equation}\label{main}\begin{array}{ll}
&\dy\sup_{\tau\in[0,\tau_1]}\int_{\mathbf{R}} \Big(\bar\r\psi^2+{\bar\r}^{\g-2}\phi^2+\phi^2_y+\psi^2_y\Big)(\tau,y)dy\\[4mm]
+&\dy\int^{\tau_1}_{0}\int_{\mathbf{R}}
\Big[\psi_y^2+{\bar\r}^{\g-2}\bar u_y\phi^2+\bar\r \bar
u_y\psi^2+\bar \r^{\g-3}\phi^2_y+\frac{\psi^2_{yy}}{\bar\r}\Big]
dyd\tau \leq C\e^{(1/2-a)}|\ln\e|^{-1/2}.
\end{array}\end{equation}
Consequently,
\begin{equation}\label{rate}
\begin{array}{ll}
\di \sup_{0\leq\tau\leq\tau_1}\|\phi(\cdot,\tau)\|_{L^\infty} \leq
\left\{\begin{array}{ll}\dy C\e^{1/6}|\ln\e|^{-1/4},~~~ &{\rm if}~~~1<\g\leq2,\\
\dy C\e^\frac{1}{\g+4}|\ln\e|^{(1-\g)/4},~~~&{\rm if}~~~ \g>2,
\end{array}\right. \\[5mm]
\di \sup_{0\leq\tau\leq\tau_1}\|\psi(\cdot,\tau)\|_{L^\infty} \leq
\left\{\begin{array}{ll}\dy C\e^{1/8}|\ln\e|^{-1/2},~~&{\rm if}~~~1<\g\leq2,\\
\dy C\e^\frac{\g+1}{4(\g+4)}|\ln\e|^{-1/2},~~&{\rm if}~~~ \g>2.
\end{array}\right.
\end{array}
\end{equation}
\end{lemma}

\

\textbf{ Proof of Lemma \ref{len}:}\ \ The proof of Lemma \ref{len}
consists of the following steps.

\underline{\it Step 1. }\quad First, define
\begin{equation*}
E:=\Phi(\r, \bar\r)+\frac{\psi^2}{2},
\end{equation*}
where
\begin{equation}\label{qua}
\Phi(\r, \bar\r):=\int^{\r}_{\bar\r}\frac{p(\xi)-p(
\bar\r)}{\xi^2}d\xi=\frac{1}{(\gamma-1)\r}
\big(p(\r)-p(\bar\r)-p^\prime(\bar\r)\phi\big).
\end{equation}
Direct computations yield
$$\begin{array}{ll}
&\Big(\r E\Big)_\tau+\Big[\r u E-\psi_y\psi+
\big(p(\r)-p(\bar\r)\big)\psi\Big]_y\\[3mm]
&+\psi_y^2 +\bar u_y\Big(p(\r)-p(\bar\r)-p'(\bar\r)
\phi\Big)+\psi^2\r \bar u_y = \bar u_{yy}\psi.
\end{array}$$
Then integrating the above equation over ${\mathbf{R}}^1\times
[0,\tau]$ and using \eqref{rup}, \eqref{eqr} and (\ref{qua}) imply
\begin{equation}\label{inten}\begin{array}{ll}
&\dy\int_{\mathbf{R}}
\Big(\bar\r\psi^2+{\bar\r}^{\g-2}\phi^2\Big)dy+
\int^\tau_{0}\int_{\mathbf{R}} \Big(\psi_y^2 +{\bar\r}^{\g-2}\bar
u_y\phi^2+\bar\r \bar u_y\psi^2\Big)dyd\tau \leq C\dy
\int^\tau_{0}\int_{\mathbf{R}} |\bar u_{yy}\psi| dyd\tau.
\end{array}\end{equation}
By Sobolev inequality and Lemma \ref{appu}, one has
\begin{equation}\label{uyy}\begin{array}{ll}
&\dy\int^\tau_{0}\int_{\mathbf{R}} |\bar u_{yy}\psi| dyd\tau\leq
C\int^\tau_{0}
\|\bar u_{yy}\|_{L^1}\|\psi\|^{1/2} \|\psi_y\|^{1/2}d\tau\\[4mm]
\leq & \dy C\int^\tau_{0}
\frac{1}{\tau+\delta/\epsilon}\|\psi\|^{1/2}
\|\psi_y\|^{1/2}d\tau\\
\leq &\dy \frac{1}{8}\int^\tau_{0} \|\psi_y\|^{2}
d\tau+C\int^\tau_{0}
(\frac{1}{\tau+\delta/\epsilon})^{4/3}\|\psi\|^{2/3}d\tau\\[4mm]
\leq & \dy \frac{1}{8}\int^\tau_{0}\|\psi_y\|^{2}
d\tau+\frac{1}{8}\sup_{\tau \in[0,\tau_1]}\|\sqrt{\bar\r}\psi\|^{2}+
C\Big(\mu^{-1/3}\int^\infty_{0}(\frac{1}{\tau+\delta/\epsilon})^{4/3}
d\tau\Big)^{3/2}\\[4mm]\leq & \dy
\frac{1}{8}\int^\tau_{0}\|\psi_y\|^{2}
d\tau+\frac{1}{8}\sup_{\tau\in[0,\tau_1]}\|\sqrt{\bar\r}\psi\|^{2}+
C\big(\f{\e}{\mu\d}\big)^{1/2}.
\end{array}\end{equation}
Combining (\ref{inten}) and (\ref{uyy}) and recalling \eqref{mu}
yield
\begin{equation}\label{energ}
\begin{array}{ll}
\dy\sup_{\tau\in[0,\tau_1]}\int_{\mathbf{R}}
\Big(\bar\r\psi^2+{\bar\r}^{\g-2}\phi^2\Big)(\tau,y)dy+
\int^{\tau_1}_{0}\int_{\mathbf{R}} \Big[\psi_y^2+{\bar\r}^{\g-2}\bar
u_y\phi^2+\bar\r \bar u_y\psi^2\Big] dyd\tau\\
\dy\hspace{9cm} \leq C\e^{(1/2-a)}|\ln\e|^{-1/2}.
\end{array}\end{equation}

\underline{\it Step 2. }\quad We make use of the idea in \cite{hq}
with modifications to derive the estimation of $\phi_{y}$.
Differentiating $(\ref {mass})$ with respect to $y$ and then
multiplying the resulted equation by $\phi_y/\r^3$ to get
\begin{equation}\label{m}\begin{array}{ll}
&\dy(\frac{\phi_y^2}{2\r^3})_\tau +(\frac{u\phi_y^2}{2\r^3})_y
+\frac{\psi_{yy}\phi_y}{\r^2}= -\frac{\phi_y}{\r^3}({\bar
u}_{yy}\phi+{\bar\r}_{yy}\psi+2{\bar\r}_y\psi_y).
\end{array}\end{equation}
Multiplying $(\ref{mome})$ by $\phi_y/\r^2$ gives
\begin{equation}\label{v}\begin{array}{ll}
&\dy (\frac{\psi\phi_y}{\r})_\tau-(\frac{\psi\phi_\tau}{\r}
 +{\bar\r}_y\frac{\psi^2}{\r})_y-\psi_y^2+p'(\r)\frac{\phi_y^2}{\r^2}
-\frac{\psi_{yy}\phi_y}{\r^2}\\[2mm]
 &\dy-{\bar u}_y\frac{\psi_y\phi}{\r}+2{\bar\r}_y\frac{\psi\psi_y}{\r}+{\bar\r}_{yy}
 \frac{\psi^2}{\r}
 +{\bar\r}_y{\bar u}_y\frac{\psi\phi}{\r^2}-\bar\r{\bar u}_{y}
 \frac{\psi\phi_y}{\r^2}=-g\frac{\phi_y}{\r^2}.
 \end{array}\end{equation}
 Adding $(\ref m)$ and $(\ref v)$ together,
then integrating the resulted equation over ${\mathbf{R}}^1\times
[0,\tau]$ imply
\begin{equation}\label{int22}\begin{array}{ll}
&\dy\int_{\mathbf{R}}\Big(\frac{\phi_y^2}{2\r^3}+
\frac{\psi\phi_y}{\r}\Big)dy +\int^\tau_{0}\int_{\mathbf{R}}
p'(\r)\frac{\phi_y^2}{\r^2}dyd\tau
\\[5mm]
= &\dy\int^\tau_{0}\int_{\mathbf{R}}\Big\{ \psi_y^2+{\bar
u}_y\frac{\psi_y\phi}{\r}-2{\bar\r}_y\frac{\psi\psi_y}{\r}-{\bar\r}_{yy}
 \frac{\psi^2}{\r}
 -{\bar\r}_y{\bar u}_y\frac{\psi\phi}{\r^2}+\bar\r{\bar u}_{y}
 \frac{\psi\phi_y}{\r^2}\\[5mm]
 &\dy
-\frac{\phi_y}{\r^3}\Big({\bar
u}_{yy}\phi+{\bar\r}_{yy}\psi+2{\bar\r}_y\psi_y\Big)
-g\frac{\phi_y}{\r^2}\Big\}dyd\tau.
 \end{array}\end{equation}
The combination of (\ref{energ}) and (\ref{int22}) leads to
 \begin{equation}\label{int3}
 \begin{array}{ll}
&\dy\int_{\mathbf{R}}\Big(\frac{\phi_y^2}{\bar\r^3}
+\bar\r\psi^2+{\bar\r}^{\g-2}\phi^2\Big)dy+\int^\tau_{0}\int_{\mathbf{R}}
\Big(\psi_y^2+{\bar\r}^{\g-2}\bar u_y\phi^2+\bar\r \bar
u_y\psi^2+\bar\r^{\g-3}\phi_y^2\Big)dyd\tau
\\[3mm]
&\dy\leq C\int^\tau_{0}\int_{\mathbf{R}}\Big\{ |{\bar
u}_y\frac{\psi_y\phi}{\bar\r}|+|{\bar\r}_y\frac{\psi\psi_y}{\bar\r}|+|{\bar\r}_{yy}
 \frac{\psi^2}{\bar\r}|+|{\bar\r}_y{\bar u}_y\frac{\psi\phi}{\bar\r^2}|+|{\bar u}_{y}
 \frac{\psi\phi_y}{\bar\r}|+|{\bar\r}_y\frac{\phi_y}{\bar\r^3}\psi_y|\\[3mm]
 &\dy\qquad\qquad\qquad+|\frac{{\bar
 u}_{yy}}{\bar\r^3}\phi_y\phi|+|\frac{{\bar\r}_{yy}}{\bar\r^3}\phi_y\psi|+|g\frac{\phi_y}{\bar\r^2}|\Big\}dyd\tau+ C\e^{(1/2-a)}|\ln\e|^{-1/2}\\
 &\dy :=\sum_{i=1}^9I_i+C\e^{(1/2-a)}|\ln\e|^{-1/2}.
 \end{array}\end{equation}
Now we estimate the terms on the right hand side of (\ref{int3}) one
by one.
 By Lemma \ref{appu}, (\ref{rup}) and Cauchy's inequality, it holds
 that
\begin{equation}\label{ry1}\begin{array}{ll}
I_1&=\dy\int^\tau_{0}\int_{\mathbf{R}} |{\bar
u}_y\frac{\psi_y\phi}{\bar\r}|dyd\tau\dy\leq
\frac{1}{8}\int^\tau_{0}\|\psi_{y}\|^{2} d\tau
+C\int^\tau_{0}\int_{\mathbf{R}} {\bar\r}^{\g-2}{\bar u}_y\phi^2\frac{{\bar u}_y}{{\bar\r}^{\g}}dyd\tau\\[3mm]
&\dy\leq\frac{1}{8}\int^\tau_{0}\|\psi_{y}\|^{2} d\tau+
C\mu^{-\g}\delta^{-1}\e\int^\tau_{0}\int_{\mathbf{R}}{\bar\r}^{\g-2}
{\bar
u}_y\phi^2dyd\tau\\[2mm]
&\dy\leq\frac{1}{8}\int^\tau_{0}\|\psi_{y}\|^{2} d\tau+ \f18
\int^\tau_{0}\int_{\mathbf{R}}{\bar\r}^{\g-2} {\bar
u}_y\phi^2dyd\tau
\end{array}\end{equation}
where we have used the fact that
\begin{equation}\label{fact1}
C\mu^{-\g}\delta^{-1}\e=C\e^{1-a-a\gamma}|\ln\e|^{-\gamma} \leq
C\e^{\f12}|\ln\e|^{-\gamma}\leq \f18,\qquad {\rm if}~~ \e\ll1.
\end{equation}
From Lemma \ref{appu} (i), one has
\begin{equation}
{\bar\r}_y={\bar\r}^{\f{3-\gamma}2}{\bar u}_{y},\label{de}
\end{equation}
one can get
\begin{equation}\label{ry11}\begin{array}{ll}
I_2&\dy=\int^\tau_{0}\int_{\mathbf{R}}
|{\bar\r}_y\frac{\psi\psi_y}{\bar\r}|dyd\tau\\[3mm]
&\dy\leq \frac{1}{8}\int^\tau_{0}\|\psi_{y}\|^{2} d\tau
+C\int^\tau_{0}\int_{\mathbf{R}} {\bar\r}^{-\g}{\bar u}_y(\bar\r{\bar u}_y\psi^2)dyd\tau\\[3mm]
&\dy\leq\frac{1}{8}\int^\tau_{0}\|\psi_{y}\|^{2} d\tau+
C\mu^{-\gamma}\delta^{-1}\e\int^\tau_{0}\int_{\mathbf{R}} {\bar\r}{\bar
u}_y\psi^2dyd\tau\\
&\dy \leq\frac{1}{8}\int^\tau_{0}\|\psi_{y}\|^{2} d\tau+
\frac{1}{8}\int^\tau_{0}\int_{\mathbf{R}} {\bar\r}{\bar
u}_y\psi^2dyd\tau
\end{array}\end{equation}
where in the last inequality we used the fact \eqref{fact1}.

Recalling \eqref{der22} from Lemma \ref{appr} and the fact $(i)$ in
Lemma \ref{appu}, one can arrive at
\begin{align}\label{rr2}
\begin{array}{ll}
|{\bar\r}_{xx}|&\dy=|\frac{2}{\gamma+1}{\bar\r}^{\f{3-\gamma}2}\omega^r_{\delta
xx}+\frac{2(3-\gamma)}{(\gamma+1)^2}
{\bar\r}^{2-\gamma}(\omega^r_{\delta x})^2|\\
&\dy \leq C\Big({\bar\r}^{\f{3-\gamma}2}\frac{{\bar u}_{
x}}{\delta}+{\bar\r}^{2-\gamma}{\bar u}_{ x}^2\Big).
\end{array}\end{align}
Thus one has
\begin{equation}\label{ry2}\begin{array}{ll}
I_3&\dy=\int^\tau_{0}\int_{\mathbf{R}} |{\bar\r}_{yy}
\frac{\psi^2}{\bar\r}| dyd\tau\\[3mm]
&\leq\dy C\e\d^{-1}\int^\tau_{0}\int_{\mathbf{R}} |\bar\r{\bar
u}_y\psi^2{\bar\r}^{-\f{\g+1}{2}}|dyd\tau
+C\e\d^{-1}\int^\tau_{0}\int_{\mathbf{R}} |{\bar\r}{\bar u}_y\psi^2{\bar \r}^{-\g}|dyd\tau\\[3mm]
&\leq\dy C\mu^{-\gamma}\delta^{-1}\e\int^\tau_{0}\int_{\mathbf{R}}
{\bar\r}{\bar u}_y\psi^2dyd\tau\\
& \leq\dy \frac{1}{8}\int^\tau_{0}\int_{\mathbf{R}} {\bar\r}{\bar
u}_y\psi^2dyd\tau,\qquad{\rm if}~~\e\ll1.
\end{array}\end{equation}
It follows from Lemma $\ref{appu}$ and (\ref{de}) that
\begin{equation}\label{uy51}\begin{array}{ll}
I_4&=\dy\int^\tau_{0}\int_{\mathbf{R}} |{\bar\r}_y{\bar u}_y\frac{\psi\phi}{\bar\r^2}|dyd\tau\\[3mm]
&\dy\leq \frac{1}{8}\int^\tau_{0}\int_{\mathbf{R}}\bar\r{\bar
u}_{y}\psi^2 dy d\tau+C\int^\tau_{0}\int_{\mathbf{R}} {\bar
u}_y{\bar\r}^{\g-2}
\phi^2\f{{\bar \r}^2_y}{{\bar\r}^{3+\g}}dyd\tau\\[4mm]
&\dy\leq \frac{1}{8}\int^\tau_{0}\int_{\mathbf{R}}\bar\r{\bar
u}_{y}\psi^2 dy d\tau
+C\f{\e^2}{\mu^{2\g}\delta^{2}}\int^\tau_{0}\int_{\mathbf{R}}{\bar\r}^{\g-2} {\bar u}_y\phi^2dyd\tau\\[4mm]
&\dy\leq \frac{1}{8}\int^\tau_{0}\int_{\mathbf{R}}\bar\r{\bar
u}_{y}\psi^2 dy d\tau +\frac{1}{8}\int^\tau_{0}\int_{\mathbf{R}}
{\bar\r}^{\g-2}{\bar u}_y\phi^2dyd\tau.
\end{array}\end{equation}
Similarly, it holds that
\begin{equation}\label{uy52}\begin{array}{ll}
I_5&\dy=\int^\tau_{0}\int_{\mathbf{R}} |{\bar u}_{y}
 \frac{\psi\phi_y}{\bar\r}|dyd\tau\\[3mm]
&\dy
\leq\dy\frac{1}{8}\int^\tau_{0}\|\bar\r^{\f{\gamma-3}{2}}\phi_{y}\|^{2}
d\tau+ C\int^\tau_{0}\int_{\mathbf{R}} {\bar\r}{\bar u}_y\psi^2\f{{\bar u}_y}{{\bar\r}^{\g}}dyd\tau\\[3mm]
&\leq\dy
\frac{1}{8}\int^\tau_{0}\|\bar\r^{\f{\gamma-3}{2}}\phi_{y}\|^{2}
d\tau+ C\f{\e}{\d\mu^\g}\int^\tau_{0}\int_{\mathbf{R}} {\bar\r}{\bar u}_y\psi^2dyd\tau\\[3mm]
&\leq\dy
\frac{1}{8}\int^\tau_{0}\|\bar\r^{\f{\gamma-3}{2}}\phi_{y}\|^{2}
d\tau+ \frac{1}{8}\int^\tau_{0}\int_{\mathbf{R}} {\bar\r}{\bar
u}_y\psi^2dyd\tau.
\end{array}\end{equation}
By Lemma $\ref{appu}$, the equality (\ref{de}) and Cauchy's
inequality, one has
\begin{equation}\label{uy5}\begin{array}{ll}
I_6&\dy=\int^\tau_{0}\int_{\mathbf{R}} |{\bar\r}_y\frac{\phi_y}{\bar\r^3}\psi_y|dyd\tau\\[3mm]
&\leq\dy\frac{1}{8}\int^\tau_{0}\|\bar\r^{\f{\gamma-3}{2}}\phi_{y}\|^{2}
d\tau+C\int^\tau_{0}\int_{\mathbf{R}} \f{{\bar u}^2_y}{{\bar\r}^{2\g}}\psi_y^2dyd\tau\\[4mm]
&\leq\dy\frac{1}{8}\int^\tau_{0}\|\r^{(\gamma-3)/2}\phi_{y}\|^{2}
d\tau+C\f{\e^2}{\d^2\mu^{2\g}}\int^\tau_{0}\int_{\mathbf{R}}\psi_y^2dyd\tau\\[4mm]
&\leq\dy  \frac{1}{8}\int^\tau_{0}\|\r^{(\gamma-3)/2}\phi_{y}\|^{2}
d\tau+ \f18\int^\tau_{0}\int_{\mathbf{R}}\psi_y^2dyd\tau.
\end{array}\end{equation}
Similarly, $I_7$ can be estimated as
\begin{equation}\label{i7}\begin{array}{ll}
I_7&\dy=\int^\tau_{0}\int_{\mathbf{R}} |\frac{{\bar
 u}_{yy}}{\bar\r^3}\phi_y\phi|dyd\tau\\[3mm]
&\leq\dy\frac{1}{8}\int^\tau_{0}\|\bar\r^{\f{\gamma-3}{2}}\phi_{y}\|^{2}
d\tau+C\f{\e}{\d}\int^\tau_{0}\int_{\mathbf{R}} {\bar\r}^{\g-2}{\bar
u}_y
\phi^2|{\bar u}_{yy}|{\bar\r}^{-1-2\g}dyd\tau\\[4mm]
&\leq\dy\frac{1}{8}\int^\tau_{0}\|\bar\r^{\f{\gamma-3}{2}}\phi_{y}\|^{2}
d\tau+C\f{\e^3}{\d^3\mu^{1+2\g}}\int^\tau_{0}\int_{\mathbf{R}} {\bar
u}_y{\bar\r}^{\g-2}
\phi^2dyd\tau\\[4mm]
&\leq\dy
\frac{1}{8}\int^\tau_{0}\|\bar\r^{\f{\gamma-3}{2}}\phi_{y}\|^{2}
d\tau +\f18\int^\tau_{0}\int_{\mathbf{R}} {\bar
u}_y{\bar\r}^{\g-2} \phi^2dyd\tau.
\end{array}\end{equation}
It follows from \eqref{rr2} that
\begin{equation}
\label{uy55}
\begin{array}{ll}
I_8&=\dy\int^\tau_{0}\int_{\mathbf{R}}
|\f{{\bar\r}_{yy}}{\bar\r^3}\phi_y\psi| dy d\tau\\
&\dy \leq C\e\d^{-1}\int_0^\tau\int_{\mathbf{R}}|\bar\r^{-\f{3+\g}{2}}\bar u_y \phi_y\psi| dyd\tau+C\e\d^{-1}\int_0^\tau\int_{\mathbf{R}}|\bar\r^{-1-\g}\bar u_y \phi_y\psi| dyd\tau\\
&\leq\dy\frac{1}{8}\int^\tau_{0}\|\bar\r^{\f{\gamma-3}{2}}\phi_{y}\|^{2}
d\tau+C\e^2\d^{-2}\int^\tau_{0}\int_{\mathbf{R}}
(\bar\r^{-2\g}+\bar\r^{1-3\g})|\bar u_y|^2\psi^2dyd\tau\\[3mm]
&\leq\dy\frac{1}{8}\int^\tau_{0}\|\bar\r^{\f{\gamma-3}{2}}\phi_{y}\|^{2}
d\tau+C\e^3\d^{-3}\mu^{-3\g}\int^\tau_{0}\int_{\mathbf{R}}
\bar\r\bar u_y\psi^2dyd\tau\\
&\leq\dy\frac{1}{8}\int^\tau_{0}\|\bar\r^{\f{\gamma-3}{2}}\phi_{y}\|^{2}
d\tau+\f18\int^\tau_{0}\int_{\mathbf{R}} \bar\r\bar
u_y\psi^2dyd\tau,\qquad{\rm if}~~\e\ll1.
\end{array}
\end{equation}
Finally, one has
\begin{equation}\label{I9}
I_9=\dy\int^\tau_{0}\int_{\mathbf{R}}
|g\frac{\phi_{y}}{\bar\r^2}|dyd\tau
\leq\frac{1}{8}\int^\tau_{0}\|\bar\r^{\f{\gamma-3}{2}}\phi_{y}\|^{2}
d\tau+C\int^\tau_{0}\int_{\mathbf{R}}\f{g^2}{\bar\r^{1+\g}}dyd\tau.
\end{equation}
Recalling that (\ref{g}), (\ref{rup}) and (\ref{de}), one can get
\begin{equation}
\begin{array}{ll}
\di |g|&\di \leq|\bar
u_{yy}|+|\bar\r\bar u_y\psi |+C|\bar\r^{\g-2}\bar\r_y\phi|\\
&\di\leq |\bar u_{yy}|+|\bar\r\bar u_y\psi
|+C|\bar\r^{\f{\g-1}2}\bar u_y\phi|.
\end{array}
\end{equation}
Thus the last term in (\ref{I9}) can be estimated by
\begin{equation}\label{ry3}\begin{array}{ll}
&\dy\Big|\int^\tau_{0}\int_{\mathbf{R}}\f{g^2}{\bar\r^{1+\g}}dyd\tau\Big|\\[3mm]
\leq&\dy C\int^\tau_{0}\int_{\mathbf{R}} \f{{\bar
u}^2_{yy}}{{\bar\r}^{1+\g}}dyd\tau+C\int^\tau_{0}\int_{\mathbf{R}}
\bar\r^{1-\g}{\bar u}^2_y\psi^2dyd\tau + \dy
C\int^\tau_{0}\int_{\mathbf{R}} {\bar\r}^{-2}{\bar u}^2_y\phi^2dyd\tau\\[3mm]
\leq&\dy C\f{\e^3}{\mu^{1+\g}}\int_0^\tau\|\bar
u_{xx}\|_{L^2(dx)}^2d\tau+C\f{\e}{\d\mu^{\g}}\int^\tau_{0}\int_{\mathbf{R}}\big(\bar\r
{\bar u}_y\psi^2+
{\bar\r}^{\g-2}{\bar u}_y\phi^2\big)dyd\tau\\[3mm]
\leq&\dy C\f{\e}{\d\mu^{1+\gamma}}\int_0^\tau
(\tau+\f{\d}{\e})^{-2}d\tau+C\f{\e}{\d\mu^\gamma}\int^\tau_{0}\int_{\mathbf{R}}
\big(\bar\r{\bar u}_y\psi^2+{\bar
u}_y{\bar\r}^{\g-2}\phi^2\big)dyd\tau\\[3mm]
\leq&\dy
C\f{\e^2}{\d^2\mu^{1+\g}}+\f18\int^\tau_{0}\int_{\mathbf{R}}
\big(\bar\r{\bar u}_y\psi^2+{\bar
u}_y{\bar\r}^{\g-2}\phi^2\big)dyd\tau, \qquad{\rm if}~~\e\ll1.
\end{array}
\end{equation}
Substituting \eqref{ry1}-(\ref{ry3}) into (\ref{int3}) and recalling
\eqref{mu}, it holds that
\begin{equation}\label{deri2}\begin{array}{ll}
\dy\int_{\mathbf{R}}\Big(\bar\r\psi^2+{\bar\r}^{\g-2}\phi^2+\phi_y^2\Big)dy+
\int^\tau_{0}\int_{\mathbf{R}}\Big(\psi_{y}^2+
 {\bar\r}^{\g-2}\bar u_y\phi^2+\bar\r \bar u_y\psi^2+\bar\r^{\gamma-3}\phi_y^2\Big) dyd\tau\\[4mm]
\leq\dy C\e^{(1/2-a)}|\ln\e|^{-1/2}, \qquad{\rm if}~~\e\ll1.
\end{array}\end{equation}

\underline{\it Step 3. }\quad As the last step, we estimate
$\dy\sup_{\tau\in[0,\tau_1]}\|\psi_y\|$. For this, multiplying
$(\ref{mome})$ by $-\psi_{yy}/\r$ gives
\begin{equation}\label{vp}
 (\frac{\psi^2_y}{2})_\tau-(\psi_y\psi_\tau
 +u\frac{\psi^2_y}{2})_y+u_y\frac{\psi_y^2}{2}
 -p'(\r)\frac{\phi_y\psi_{yy}}{\r}
+\frac{\psi_{yy}^2}{\r}=g\frac{\psi_{yy}}{\r}.
\end{equation}
Integrating the above equation over ${\mathbf{R}}^1\times [0,\tau]$
yields
\begin{equation}\label{intmv}\begin{array}{ll}
&\dy\int_{\mathbf{R}}  \frac{\psi^2_y}{2}dy
+\int^\tau_{0}\int_{\mathbf{R}}\Big (\frac{{\bar u}_y\psi_y^2}{2}
  +\frac{\psi_{yy}^2}{\r}\Big)dyd\tau
=\int^\tau_{0}\int_{\mathbf{R}}\Big\{
p'(\r)\frac{{\psi}_{yy}\phi_y}{\r}
+g\frac{\psi_{yy}}{\r}-\f{\psi_y^3}{2}\Big\}dyd\tau.
 \end{array}\end{equation}
First, one has
\begin{equation}\label{uy3}\begin{array}{ll}
&\dy\Big|\int^\tau_{0}\int_{\mathbf{R}}
p'(\r)\frac{{\psi}_{yy}\phi_y}{\r}dyd\tau\Big|\leq
\frac{1}{8}\int^\tau_{0}\int_{\mathbf{R}}\frac{\psi_{yy}^{2}}{\r}
dyd\tau+C\int^\tau_{0} \|\bar\r^{\f{\gamma-3}2}\phi_y\|^{2}d\tau.
\end{array}\end{equation}
Then it follows from (\ref{ry3}) and (\ref{deri2}) that
\begin{equation}\label{uy6}
\begin{array}{rl}
\dy\Big|\int^\tau_{0}\int_{\mathbf{R}}
g\frac{\psi_{yy}}{\r}dyd\tau\Big|\leq &\dy
\frac{1}{8}\int^\tau_{0}\int_{\mathbf{R}}
\frac{\psi^2_{yy}}{\r}dyd\tau+
C\Big|\int^\tau_{0}\int_{\mathbf{R}} \frac{g^2}{\bar\r}dyd\tau\Big|\\[3mm]
\leq&\dy \frac{1}{8}\int^\tau_{0}\int_{\mathbf{R}}
\frac{\psi^2_{yy}}{\r}dyd\tau+ C\e^{1/2-a}|\ln\e|^{-1/2}.
\end{array}
\end{equation}
Furthermore, we can compute that
\begin{equation}\label{uy7}
\begin{array}{ll}
\dy\Big|\int^\tau_{0}\int_{\mathbf{R}}
\frac{\psi_y^3}{2}dyd\tau\Big|&\dy \leq
C\int^\tau_{0}\|\psi_{yy}\|^{\f12}\|\psi_y\|^{\f52} d\tau\\
&\dy \leq
\frac{1}{8}\int^\tau_{0}\int_{\mathbf{R}}\frac{\psi_{yy}^{2}}{\r}
dyd\tau+C\int^\tau_{0} \|\psi_y\|^{\f{10}{3}}d\tau\\
&\dy \leq
\frac{1}{8}\int^\tau_{0}\int_{\mathbf{R}}\frac{\psi_{yy}^{2}}{\r}
dyd\tau+C\sup_{\tau\in[0,\tau_1]}\|\psi_y\|^{\f43}\int^\tau_{0}
\|\psi_y\|^2d\tau\\
&\dy \leq
\frac{1}{8}\int^\tau_{0}\int_{\mathbf{R}}\frac{\psi_{yy}^{2}}{\r}
dyd\tau+C\int^\tau_{0} \|\psi_y\|^2d\tau,
\end{array}
\end{equation}
where in the last inequality we have used the a priori assumption
\begin{equation}\label{assump-1}
\sup_{\tau\in[0,\tau_1]}\|\psi_y\|\leq 1.
\end{equation}
Substituting \eqref{uy3}, \eqref{uy6} and (\ref{uy7}) into
(\ref{intmv}) and using (\ref{rup}) and (\ref{deri2}), it holds that
\begin{equation}\label{deri1}
\begin{array}{ll}
&\dy\int_{\mathbf{R}}\psi_y^2dy+ \int^\tau_{0}\int_{\mathbf{R}}\Big
({\bar u}_y\psi_y^2
  +\frac{\psi_{yy}^2}{\bar\r}\Big)dyd\tau\leq\dy C\e^{(1/2-a)}|\ln\e|^{-1/2}.
\end{array}\end{equation}
Therefore, \eqref{main} can be derived directly from (\ref{deri2})
and (\ref{deri1}) and the a priori assumption \eqref{assump-1} is
verified if $\e$ is suitably small. It follows from \eqref{main}
that if $1<\gamma\leq2$, then
\begin{equation}
\label{a1}
\begin{array}{ll}
\dy\sup_{0\leq\tau\leq\tau_1}\|\phi(\cdot,\tau)\|_{L^\infty}&\dy \leq \sqrt2\sup_{0\leq\tau\leq\tau_1}\|\phi(\cdot,\tau)\|^{1/2}\|\phi_y(\cdot,\tau)\|^{1/2}\\
&\dy\leq
C\sup_{0\leq\tau\leq\tau_1}\big(\int_{\mathbf{R}}\bar\r^{\g-2}\phi^2dy\big)^{\f14}
\big(\int_{\mathbf{R}}\phi_y^2dy\big)^{\f14}\\
&\dy \leq C\e^{1/6}|\ln\e|^{-1/4},
\end{array}
\end{equation}
and if $\gamma>2$, then
\begin{equation}\label{a2}
\begin{array}{ll}
\dy\sup_{0\leq\tau\leq\tau_1}\|\phi(\cdot,\tau)\|_{L^\infty}&\dy \leq \sqrt2\sup_{0\leq\tau\leq\tau_1}\|\phi(\cdot,\tau)\|^{1/2}\|\phi_y(\cdot,\tau)\|^{1/2}\\
&\dy \leq
C\sup_{0\leq\tau\leq\tau_1}\big(\int_{\mathbf{R}}\mu^{2-\g}\bar\r^{\g-2}\phi^2dy\big)^{\f14}
\big(\int_{\mathbf{R}}\phi_y^2dy\big)^{\f14}
\\
&\dy\leq C\mu^{\f12-\f\gamma4}\e^{\f14-\frac{1}{2(\g+4)}}|\ln\e|^{-\f14}\\
&\dy\leq C\e^{\frac{1}{\g+4}}|\ln\e|^{\frac{1-\gamma}{4}}.
\end{array}
\end{equation}
And we also have
\begin{equation}\label{a3}
\sup_{\tau\in[0,\tau_1]}\|\psi_y\| \leq C\e^{(\f14-\frac{a}{2})}|\ln\e|^{-1/4}\leq 1.
\end{equation}
 So from \eqref{a1}-(\ref{a3}) the a priori
assumption \eqref{assump} is verified if $\e\ll 1$. On the other
hand, by using Sobolev inequality, one can get
\begin{equation}
\begin{array}{ll}
\dy\sup_{0\leq\tau\leq\tau_1}\|\psi(\cdot,\tau)\|_{L^\infty} &\di
\leq\sqrt2\sup_{0\leq\tau\leq\tau_1}\|\psi(\cdot,\tau)\|^{1/2}\|\psi_y(\cdot,\tau)\|^{1/2}\\[3mm]
&\di \leq
\sqrt2\mu^{-\f14}\sup_{0\leq\tau\leq\tau_1}\|\sqrt{\bar\r}\psi(\cdot,\tau)\|^{1/2}\|\psi_y(\cdot,\tau)\|^{1/2}\\
&\di \leq
\left\{\begin{array}{ll}\dy C\e^{1/8}|\ln\e|^{-1/2},~~&{\rm if}~~~1<\g\leq2,\\
\dy C\e^\frac{\g+1}{4(\g+4)}|\ln\e|^{-1/2},~~&{\rm if}~~~ \g>2,
\end{array}\right.
\end{array}
\end{equation}
Thus the convergence rate \eqref{rate} is justified and the proof of
Lemma \ref{len} is completed.

 \hfill $\Box$
\vspace{2mm}

\textbf{ Proof of Theorem ~\ref{thm1}:}\ \ It remains to prove
\eqref{decay-rate} with $a, b$ given in \eqref{aaa} and \eqref{bbb}
respectively. From Lemma \ref{cut-off}, Lemma \ref{appu} (iii) and
Theorem \ref{thm31} and recalling that
$\mu=\e^{a}|\ln\e|,~~\delta=\e^{a}$, it holds that for any given
positive constant $h$, there exist a constant $C_h>0$ which is
independent of $\epsilon$ such that
\begin{equation*}\begin{array}{ll}
&\dy\sup_{t\geq h}\| \r(\cdot,t)-\rho^r(\frac{\cdot}{t})\|_{L^\infty}\\
\leq&\dy\sup_{0\leq\tau\leq+\i}\|
\phi(\cdot,\tau)\|_{L^\infty}+\sup_{t\geq
h}\|\bar\r(\cdot,t)-\rho^r_\mu(\frac{\cdot}{t})\|
_{L^\infty}+\sup_{t\geq h}\| \rho^r_\mu(\frac{\cdot}{t})- \rho^r(\frac{\cdot}{t})\|_{L^\infty}\\
\leq &\dy
\left\{\begin{array}{ll}\dy C_h\Big(\e^{1/6}|\ln\e|^{-1/4}+\delta|\ln\delta|+\mu\Big),~~&{\rm if}~~~1<\g\leq2,\\
\dy
C_h\Big(\e^\frac{1}{\g+4}|\ln\e|^{(1-\g)/4}+\delta|\ln\delta|+\mu\Big),~~&{\rm
if}~~~ \g>2,
\end{array}\right.\\[2mm]
\leq &\dy C_h\epsilon^{a}|\ln\epsilon |,
\end{array}\end{equation*}
and
\begin{equation*}\begin{array}{ll}
&\dy\sup_{t\geq h}\| m(\cdot,t)-m^r(\frac{\cdot}{t})\|_{L^\infty}\\
\leq&\dy \sup_{t\geq h}\Big(\|m(\cdot,t)-\bar
m(\cdot,t)\|_{L^\infty}+\|\bar m(\cdot,t)-m^r_\mu(\frac{\cdot}{t})\|
_{L^\infty}+\| m^r_\mu(\frac{\cdot}{t})- m^r(\frac{\cdot}{t})\|_{L^\infty}\Big)\\
\leq&\dy
C\sup_{0\leq\tau\leq+\i}\Big(\|\psi(\cdot,\tau)\|_{L^\infty}+\|
\phi(\cdot,\tau)\|_{L^\infty}\Big)\\
&\dy\qquad\qquad\qquad\qquad+\sup_{t\geq h}\Big(\|\bar
m(\cdot,t)-m^r_\mu(\frac{\cdot}{t})\|
_{L^\infty}+\| m^r_\mu(\frac{\cdot}{t})- m^r(\frac{\cdot}{t})\|_{L^\infty}\Big)\\
\leq &\dy
\left\{\begin{array}{ll}\dy C_h\Big(\e^{1/8}|\ln\e|^{-1/2}+\e^{1/6}|\ln\e|^{-1/4}+\delta|\ln\delta|+\mu\Big),~~&{\rm if}~~~1<\g\leq2,\\
\dy
C_h\Big(\e^\frac{\g+1}{4(\g+4)}|\ln\e|^{-1/2}+\e^\frac{1}{\g+4}|\ln\e|^{(1-\g)/4}+\delta|\ln\delta|+\mu\Big),~~&{\rm
if}~~~ \g>2,
\end{array}\right.\\[6mm]
\leq &\dy
\left\{\begin{array}{ll}\dy C_h\epsilon^{b} |\ln\epsilon |^{-\f12}, \,~~~ &{\rm if}~~~1<\g<3,\\
\dy C_h\epsilon^{\frac{1}{\g+4}} |\ln\epsilon |,~~~&{\rm if}~~~
\g\geq3,
\end{array}\right.
\end{array}\end{equation*}
 Thus
the proof of Theorem 1.1 is completed.
 \hfill $\Box$
\vspace{2mm}

\section*{Acknowledgments}
The authors would like to thank the referees for the valuable
comments and suggestions which greatly improved the presentation of
the paper. The research of F. M. Huang was supported in part by NSFC
Grant No. 10825102 for distinguished youth scholar, and National
Basic Research Program of China (973 Program) under Grant No.
2011CB808002. The research of Y. Wang was supported by the NSFC
Grant No. 10801128.

\textbf{    }
\end{CJK}
\end{document}